\documentclass[11pt]{amsart}

\usepackage{amsmath, amsfonts, amssymb,amsthm}
\usepackage{amstext}
\usepackage{mathrsfs}
\usepackage{float,epsf,subfigure}
\usepackage[all,cmtip]{xy}
\usepackage{hyperref}
\usepackage{algorithm}
\usepackage{algorithmic}
\usepackage{mathtools}
\usepackage{float}
\usepackage{bm}
\usepackage{graphicx} 
\usepackage{caption} 
\usepackage{multirow}

\usepackage[
    a4paper,
    left=3.3cm,
    right=3.3cm,
    top=3.8cm,
    bottom=3.8cm
]{geometry}
\usepackage{xcolor}
\newtheorem{thm}{Theorem}[section]

\newtheorem{def-thm}[thm]{Definition-Theorem}
\newtheorem{lem}[thm]{Lemma}

\newtheorem*{tha}{Theorem A}
\newtheorem*{thb}{Theorem B}

\theoremstyle{plain}
\theoremstyle{definition}

\def\min{\mathop{\mathrm{min}}}
\def\C{\mathbb C}

\begin{document}
\title[The periodicity of meromorphic functions]{On the periodicity of meromorphic functions of the form $ff^{(k)}$}
\author[S.S. Yang \& D.H. Zhao \& C.C. Yang \& L.W. Liao] {Shuang-Shuang Yang \& Dong-hai Zhao* \& Chung-Chun Yang \& Liang-Wen Liao}

\address{School of Mathematical Sciences \\ Suzhou University of Science and Technology \\ Suzhou, P. R. China}
\email{ssyangmath@126.com}
\address{School of Mathematics \\ Nanjing University \\ Nanjing, P. R. China}
\email{donghaizhao2021@163.com}
\address{School of Mathematics \\ Nanjing University \\ Nanjing, P. R. China}
\email{maccyang@163.com}
\address{School of Mathematics \\ Nanjing University \\ Nanjing, P. R. China}
\email{maliao@nju.edu.cn}

\thanks{This research is supported by the National Natural Science Foundation of China (Grant No. 12471077, 12601138).}
\thanks{*Corresponding author}
\date{\today}

\maketitle \thispagestyle{empty} \setcounter{page}{1}

\begin{abstract}
For a non-entire meromorphic function $f$ and an integer $k$, we will show that whenever $ff^{(k)}$ is periodic so is $f$. This result gives an affirmative answer to a significant class of the so called generalized Yang's conjecture.
\end{abstract}

\vskip\baselineskip

\setlength\arraycolsep{2pt}

\section{Introduction and main results}
We recall Yang's conjecture, whose precise statement is as follows: given a holomorphic function $f$, if $ff^{(k)}$ $(k\in\mathbb N^+)$ is periodic, is $f$ necessarily periodic? Recently, this question has attracted extensive attention from researchers. Wang and Hu \cite{wh2018} proved this conjecture is true when $k=1$. However, the case $k\geq2$ remains open as of now. Later, several scholars have started to study its differential or difference variants, as well as the case where $f$ is meromorphic, see \cite{lz2022,lly2019,lwy2020,ly2019,lk2020,lrl2022,mswy2025}. These research lead to the so called generalized Yang's conjecture, which asks whether $f$ is periodic whenever $f$ is meromorphic and $f^nf^{(k)}$ $(n,k\in\mathbb N^+)$ is periodic.  Some investigations were mentioned in \cite{mswy2025}. And main results from \cite{mswy2025} studied this question for $f^nf^{(k)}$ when $f$ is non-holomorphic meromorphic. Specifically, 
\begin{tha}[\cite{mswy2025}]
Let $n$ be a positive integer and $k$ be an even integer. Assume $f$ is transcendental meromorphic with at least one pole. If $f^nf^{(k)}$ is $T$-periodic, then $f$ is a periodic function with period $pT$, where $p\in\{1.2,\cdots,n+1\}$.
\end{tha}

\begin{thb}[\cite{mswy2025}]
Let $m,n\geq2$ be integers and $k$ be an odd integer with $k>m(n-1)$. Assume $f$ is transcendental meromorphic such that the minimum multiplicity of all its
poles is $m$. If $f^nf^{(k)}$ is $T$-periodic, then $f$ is a periodic function with period $pT$, where $p\in\{1.2,\cdots,n+1\}$.
\end{thb}

Through observation, they provided an affirmative answer to $f^nf^{(k)}$ for even $k$ and for meromorphic $f$ having at least one pole. However, in the odd case, they imposed restrictions on $n$ and the multiplicities of the poles. In the following, we consider $n= 1$ and allow poles of any positive multiplicity. For $k= 1$, the conclusion is already known for transcendental meromorphic functions, including holomorphic functions, see \cite{lrl2022}. The main result of the present paper treats arbitrary positive odd $k$. And we do not impose conditions on the growth of $f$.

Throughout this paper, a meromorphic function means a function meromorphic on $\C$, unless some other domain is specified. A $c$-period of $f$ is a number $c\neq0$ satisfying 
$$f(z+c)=f(z),$$
which need not be a smallest period.

\begin{thm}\label{thmoddk}
Let $k$ be an odd integer. Assume that $f$ is meromorphic with at least one pole. If $ff^{(k)}$ is $c$-periodic, then $f$ is a periodic function with period $c$ or $2c$. 
\end{thm}
Together with Theorem A, we immediately deduce the following conclusion.
\begin{thm}\label{thmanyk}
Assume that $f$ is meromorphic with at least one pole. If $ff^{(k)}$ is $c$-periodic, then $f$ is a periodic function with period $c$ or $2c$.
\end{thm}
Two cases in Theorem~\ref{thmanyk} do occur. To illustrate, when $f(z)=1/\sin{(\pi z/c)}$, $ff^{(k)}$ has period $c$, whereas $f$ has period $2c$ rather than $c$; when $f(z)=\cot{(\pi z/c)}$, both $ff^{(k)}$ and $f$ have period $c$.  


\section{Preliminaries}\label{sec:prelim}
This section presents several lemmas that will be needed. For convenience, we define some notations as follows. For $z_0\in\C$ and $r>0$, write
$$D(a,r)=\{z\in\C:|z-a|<r\}, \qquad D^*(a,r)=\{z\in\C:0<|z-a|<r\}$$
For a meromorphic function $u$, define
$$\mathcal{Q}_k[u]=\frac{1}{2} \sum_{l=0}^{k-1} (-1)^l u^{(l)} u^{(k-1-l)}.$$
It is not hard to check that
$$(\mathcal{Q}_k[u])'=uu^{(k)}, \qquad \mathcal{Q}_k[-u] = \mathcal{Q}_k[u].$$

A family of holomorphic functions on a domain is called \textit{normal} if every sequence has a subsequence that converges uniformly on compact subsets to a holomorphic function, or tends uniformly to infinity on compact subsets. And a family is \textit{locally bounded} if it has a common bound on each compact subset of its domain. We refer to \cite{zalcman1998normal} for more details.

\begin{lem}[\cite{zalcman1998normal}]\label{lemzalcman}
Let $\mathcal{F}$ be a family of holomorphic functions on a domain $D\in\C$. If $\mathcal{F}$ is not normal on $D$, then there exist a sequence $\{f_n\}\subset\mathcal{F}$, a sequence $\{z_n\}$ with $\lim_{n\rightarrow\infty}z_n=z_0\in D$, and a sequence of positive number $\{\varepsilon_n\}$ with $\lim_{n\rightarrow\infty}\varepsilon_n=0$ such that 
$$g_n(\xi) = f_n(z_n + \varepsilon_n\xi)$$ 
converges uniformly to a nonconstant holomorphic function $g(\xi)$ on any compact subset of $\mathbb{C}$.
\end{lem}

\begin{lem}\label{lemnormal}
Let $D\subset\C$ be a domain and $\mathcal{F}$ be a family of zero-free holomorphic functions on $D$. If
$$\big\{ff^{(k)}:f\in\mathcal{F}\big\}$$
is locally bounded on $D$, then $\mathcal{F}$ is normal on $D$.
\end{lem}
\begin{proof}
On the contrary, suppose that $\mathcal{F}$ is not normal. Applying Lemma~\ref{lemzalcman}, we can get a sequence $\{f_n\}\subset\mathcal{F}$ satisfying
$$g_n(\xi) = f_n(z_n + \varepsilon_n\xi) \,\Longrightarrow \, g(\xi)$$ 
on compact subsets of $\C$, where $\lim_{n\rightarrow\infty}z_n=z_0\in D$, $\lim_{n\rightarrow\infty}\varepsilon_n=0$ and $g$ is a nonconstant holomorphic function.

We claim that $g$ has no zeros. Suppose that $\xi_0$ is a zero of $g$. Choose $r>0$ such that $g$ has no zeros on $|\xi-\xi_0|=r$. Then, for $n$ sufficiently large, we have
$$\max_{|\xi-\xi_0|=r}|g_n(\xi)-g(\xi)|<\min_{|\xi-\xi_0|=r}|g(\xi)|.$$
Rouch\'{e}'s theorem implies that $g_n$ and $g$ have the same number of zeros inside the circle $|\xi-\xi_0|=r$ counting multiplicities. This
contradicts that every $g_n$ has no zeros. Hence $g$ is zero-free.

Cauchy’s integral formula shows that all derivatives of each order
converge uniformly on compact subsets, i.e., $g_n^{(k)}(\xi)\Longrightarrow g^{(k)}(\xi)$. Take a compact set $K\subset D$. Since $ff^{(k)}$ is locally bounded, a direct computation gives 
$$g_n(\xi)g_n^{(k)}(\xi)=\varepsilon_n^k\big(f_nf_n^{(k)}\big)(z_n + \varepsilon_n\xi)\Longrightarrow 0$$
on $K$. Furthermore, it yields
$$gg^{(k)}\equiv0$$
on $\C$. Since $g$ has no zeros, $g^{(k)}\equiv0$, which means $g$ is a polynomial of degree at most $k-1$. However, a nonconstant polynomial must have zeros. Therefore $g$ is constant, a contradiction. The assumption is not valid, i.e., $\mathcal{F}$ is normal on $D$.
\end{proof}

\begin{lem}\label{lembounded}
Assume that $f$ is a meromorphic function with a pole at $z_0$ and  $ff^{(k)}$ is $c$-periodic. Denote
$$f_j(z)=f(z+jc),\qquad j\in\mathbb N.$$
Then there exists $r>0$ such that the family $\{f_j\}$ is holomorphic, zero-free and locally bounded on $D^*(z_0,r)$. Moreover, for any $\zeta\in D^*(z_0,r)$ and any $l\in\mathbb N$, we have
$$\sup_{j\in\mathbb N}\left|f_j^{(l)}(\zeta)\right|<\infty.$$
\end{lem}
\begin{proof}
Set $F=ff^{(k)}.$ Suppose that the pole $z_0$ of $f$ has order $m$. Then, in a neighborhood of $z_0$, we have
$$f(z)= A\tau^{-m}+O\big(\tau^{-m+1}\big),$$
where $A\neq0,$ $m\in\mathbb N,$ and $\tau=z-z_0.$ Differentiation gives
$$f^{(k)}(z)= (-1)^k\frac{(m+k-1)!}{(m-1)!}A\tau^{-m-k}+O\big(\tau^{-m-k+1}\big).$$
Therefore,
\begin{equation}\label{eqFf}
F(z)=(-1)^k\frac{(m+k-1)!}{(m-1)!}A^2\tau^{-2m+k}+O\big(\tau^{-2m-k+1}\big).
\end{equation}
It follows that $z_0$ is also a pole of $F$, of order $2m+k$. Then, for any $j\in\mathbb N$, a simple calculation yields
\begin{equation}\label{eqFfj}
f_jf^{(k)}_j(z)=F(z+jc)=F(z),
\end{equation}
which indicates $z_0$ is also a pole of $f_j$. Let its order be $m_j$ and its leading coefficient of Laurent expansion of $f_j$ be $A_j$. Differentiating the Laurent expansion of $f_j$ $k$ times, substituting the expansions $f_j$ and $f_j^{(k)}$ into \eqref{eqFfj}, and then comparing the order and the leading coefficient with \eqref{eqFf} in turn, we obtain
$$m_j=m, \qquad A_j=\pm A.$$
Take $r$ so that $F$ is holomorphic and nonzero on $D^*(z_0,r)$. It yields that $f_j$ is holomorphic and zero-free on $D^*(z_0,r)$ as well. In addition, all $f_jf_j^{(k)}$ are locally bounded on $D^*(z_0,r)$ due to \eqref{eqFfj}. Employing Lemma~\ref{lemnormal}, we can see $\{f_j\}$ is normal on $D^*(z_0,r)$. Furthermore, we claim $\{f_j\}$ is local boundedness on $D^*(z_0,r)$. If not, then there exists a subsequence $\{f_{j_n}\}$  tending uniformly to infinity on a compact subset $K$ of $D^*(z_0,r)$. Denote 
$$V_j(z)=(z-z_0)^mf_j(z),$$
which is a holomorphic and zero-free function on $D(z_0,r)$, and 
$$V_j(z_0)=A_j,\qquad |V_j(z_0)|=|A|.$$
In particular, $1/V_j$ is also holomorphic on $D(z_0,r)$.  Fix $\tilde r$ with $\{z:|z-z_0|=\tilde r\}\subset K$. Then
$$\max_{|z-z_0|=\tilde r}\left|\frac{1}{f_{j_n}(z)}\right|\longrightarrow0.$$
Applying maximum modulus principle, we deduce that 
$$\frac{1}{|A|}=\left|\frac{1}{V_{j_n}(z_0)}\right|\leq\max_{|z-z_0|=\tilde r}\left|\frac{1}{V_{j_n}(z)}\right|=r^{-m}\max_{|z-z_0|=\tilde r}\left|\frac{1}{f_{j_n}(z)}\right|\longrightarrow0,\qquad n\longrightarrow\infty,$$
which is impossible. Thus, normality implies local boundedness. 

For any $\zeta\in D^*(z_0,r)$, and choose $\rho>0$ such that $\overline{D(\zeta,\rho)}\subset D^*(z_0,r)$. According to the local
boundedness, there exists a constant $M>0$ for which
$$\left|f_j(z)\right|\leq M, \qquad z\in\overline{D(\zeta,\rho)}.$$
It follows Cauchy’s integral formula that
\begin{equation}\label{eqcauchy}
\left|f_j^{(l)}(\zeta)\right|\leq \frac{l!M}{\rho^l}.    
\end{equation}
Thus, for any $l\in\mathbb N$, we have
$$\sup_{j\in\mathbb N}\left|f_j^{(l)}(\zeta)\right|<\infty.$$
\end{proof}

\begin{lem}\label{lemuv}
Let \(k\) be odd. Assume that  $f$ and $g$ are meromorphic functions with the same order $m\geq1$  and the same leading Laurent  coefficient $a \neq 0$ at pole $z_0$. If
\begin{equation}\label{eqpropcon}
ff^{(k)}=gg^{(k)}, \qquad \mathcal{Q}_k[f] = \mathcal{Q}_k[g],
\end{equation}
then $f=g$ in a neighborhood of $z_0$.
\end{lem}
\begin{proof}
Suppose, on the contrary, that $h=g-f$ is not identically zero in a neighborhood of $z_0$. Since the leading  terms of Laurent expansions $f$ and $g$ agree, there are an integer $n>-m$ and a constant $b\neq0$ such that
\begin{align*}
f(z) &= a\tau^{-m}+O(\tau^{-m+1}), \\
h(z) &= b\tau^n+O(\tau^{n+1}).    
\end{align*}    
Combining the first identity in \eqref{eqpropcon} with \(g = f + h\), we have
\begin{equation}\label{equukhhk}
0 = f h^{(k)} + h f^{(k)} + h h^{(k)}.
\end{equation}
A straightforward computation yields
\begin{align*}
f h^{(k)} &=-\frac{(m
+k-1)!}{(m-1)!} a b \tau^{n-m-k} + O(\tau^{n-m-k+1}), \\
h f^{(k)} &= n(n-1)\cdots(n-k+1)a b  \tau^{n-m-k} + O(\tau^{n-m-k+1}), \\
h h^{(k)} &= O(\tau^{2n-k}).
\end{align*}
Comparing the lowest-order terms in \eqref{equukhhk} yield $n=m+k-1\geq k$. Further calculation leads to
\begin{align*}
\mathcal{Q}_k[g]-\mathcal{Q}_k[f] & =\mathcal{Q}_k[f+h]-\mathcal{Q}_k[f] \\
& = \frac{1}{2}\sum_{l=0}^{k-1} (-1)^l \big( f^{(l)} h^{(k-1-l)} + h^{(l)} f^{(k-1-l)} \big) + \mathcal{Q}_k[h] \\
& = \sum_{l=0}^{k-1} (-1)^l f^{(l)} h^{(k-1-l)} + \mathcal{Q}_k[h] \\
& = \sum_{l=0}^{k-1}(-1)^l(-1)^l\frac{(m+l-1)!}{(m-1)!}\frac{n!}{(n-k+1+l)!}ab+O(\tau) \\
& = ab\frac{n!}{(m-1)!} \sum_{l=0}^{k-1} \frac{1}{m+l}+O(\tau).
\end{align*}
In terms of the second identity in \eqref{eqpropcon}, we conclude that
$$0=ab\frac{n!}{(m-1)!} \sum_{l=0}^{k-1} \frac{1}{m+l},$$
which is impossible. Hence $h$ must vanish identically near $z_0$. This proves $f=g$.
\end{proof}

\section{Proof of Theorem~\ref{thmoddk}}
\begin{proof}
Assume that $f$ has a pole $z_0$ of order $m$. In a neighborhood of $z_0$, we have
$$f(z)= a\tau^{-m}+O\big(\tau^{-m+1}\big),$$
where $a\neq0$ and $\tau=z-z_0.$
Lemma~\ref{lemnormal} shows that $f(z+c)$
has the same order at $z_0$ and its leading coefficient is $a_1=\pm a$. Set
$$\kappa\in\frac{a_1}{a}\in\{1,-1\},\qquad g(z)=\kappa f(z+c).$$
Then $g$ and $f$ have the same leading coefficient $a$ at $z_0$ and 
$$gg^{(k)}=ff^{(k)}.$$
The periodicity of $F=ff^{(k)}$ leads to 
$$\frac{d}{dz}\left(\mathcal{Q}_k[f](z+c)-\mathcal{Q}_k[f](z)\right)=F(z+c)-F(z)=0.$$
Hence, we obtain
\begin{equation}\label{QQC0}
\mathcal{Q}_k[f](z+c)-\mathcal{Q}_k[f](z)= C_0,\qquad C_0\in\C.
\end{equation}
Set
$$f_j(z)=f(z+jc),\qquad j\in\mathbb N.$$
Employing Lemma~\ref{lembounded}, there exists $r>0$ such that for \(\zeta \in D^*(z_0, r)\), 
$$|\mathcal{Q}_k[f_j](\zeta)|\leq \frac{1}{2} \sum_{l=0}^{k-1}\left|f_j^{(l)}(\zeta)\right| \, \left|f_j^{(k-1-l)}(\zeta)\right|\leq \frac{M^2}{2\rho^{k-1}} \sum_{l=0}^{k-1} l! (k - 1 - l)!,$$
where $M,\rho$ satisfies
$$\left|f_j(z)\right|\leq M,\qquad z\in\overline{D(\zeta,\rho)}\subset D^*(z_0,r).$$
On the other hand, \eqref{QQC0} gives
$$\mathcal{Q}_k[f_j](\zeta)=\mathcal{Q}_k[f](\zeta+jc)=\mathcal{Q}_k[f](\zeta)+jC_0.$$
If $C_0 \neq 0$, then
$$|\mathcal{Q}_k[f](\zeta)+jC_0|\geq j|C_0| - |\mathcal{Q}_k[f](\zeta)|\longrightarrow \infty,\qquad j\longrightarrow\infty,$$
a contradiction. Therefore $C_0 = 0$, which means
$$\mathcal{Q}_k[f]=\mathcal{Q}_k[f](z+c)=\mathcal{Q}_k[g](z).$$
Lemma~\ref{lemuv} shows that $g=f$ near $z_0$, i.e. $\kappa f(z + c) = f(z)$. 
If $\kappa=1$, then $c$ is a period of $f$. If $\kappa=-1$, then $2c$ is a period of $f$. This completes the proof.
\end{proof}

\bibliographystyle{plain}
\bibliography{references}

\end{document}